\documentclass[a4paper]{amsart}

\usepackage[matrix,arrow,curve,tips]{xy}
           \SelectTips{cm}{}
\usepackage{mathrsfs}

\newtheorem{dfn}{Definition}[section]
\newtheorem{prop}[dfn]{Proposition}
\newtheorem{theo}[dfn]{Theorem}
\newtheorem{cor}[dfn]{Corollary}
\newtheorem{lem}[dfn]{Lemma}
\newtheorem{rem}[dfn]{Remark}
\newtheorem{ex}[dfn]{Example}

\newcommand{\com}{\mathbin{{\scriptstyle \circ }}}

\newcommand{\RR}{\mathbb{R}}

\newcommand{\fg}{\mathfrak{g}}
\newcommand{\fp}{\mathfrak{p}}

\newcommand{\id}{\mathord{\mathrm{id}}}
\newcommand{\pr}{\mathord{\mathrm{pr}}}

\newcommand{\eC}{\mathord{\mathit{C}^{\infty}}}

\newcommand{\Hom}{\mathord{\mathrm{Hom}}}

\newcommand{\Der}{\mathord{\mathrm{Der}}}
\newcommand{\End}{\mathord{\mathrm{End}}}

\newcommand{\Prim}{\mathord{\mathscr{P}}}
\newcommand{\Univ}{\mathord{\mathscr{U}}}
\newcommand{\gr}{\mathord{\mathrm{gr}}}

\newcommand{\ten}{\mathbin{\otimes}}
\newcommand{\cu}{\mathord{\epsilon}}
\newcommand{\cm}{\mathord{\Delta}}

\newcommand{\LRA}{\mathsf{LieRinAlg}}
\newcommand{\RBA}{\mathsf{RinBiAlg}}

\begin{document}

\title{On the universal enveloping algebra of a Lie-Rinehart algebra}

\author{I. Moerdijk}
\address{Mathematical Institute, Utrecht University,
         P.O. Box 80.010, 3508 TA Utrecht, The Netherlands}
\email{moerdijk@math.uu.nl}

\author{J. Mr\v{c}un}
\address{Department of Mathematics, University of Ljubljana,
         Jadranska 19, 1000 Ljubljana, Slovenia}
\email{janez.mrcun@fmf.uni-lj.si}

\thanks{The second author was supported in part by
        the Slovenian Ministry of Science}


\begin{abstract}
We review the extent to which the universal enveloping algebra
of a Lie-Rinehart algebra resembles a Hopf algebra, and
refer to this structure as a Rinehart bialgebra. We then prove
a Cartier-Milnor-Moore type theorem for such Rinehart
bialgebras.
\end{abstract}

\maketitle

\section{Lie-Rinehart algebras} \label{sec1}

Throughout this paper $k$ will denote a field
of characteristic $0$.
In this section we will briefly recall the definition
of Lie-Rinehart algebra and mention some basic examples.
Lie-Rinehart algebras were introduced and studied
by Herz \cite{Her}, Palais \cite{Pal} and Rinehart \cite{Rin}.
They form the algebraic counterpart of the more geometric
notion of Lie algebroid, which has become better known.
It was Huebschmann \cite{Hue} who gave Lie-Rinehart algebras their name,
and also emphasized their advantages over Lie algebroids
in contexts where singularities arise.
See also \cite{Hue2} for an excellent survey.

Let $R$ be a unital {\em commutative} algebra over $k$.
A {\em Lie-Rinehart algebra over} $R$ is a Lie algebra
$L$ over $k$, equipped with a structure of a unital left $R$-module
and a homomorphism of Lie algebras $\rho\!:L\to\Der_{k}(R)$
into the Lie algebra of derivations on $R$,
which is a map of left $R$-modules
satisfying the
following Leibniz rule
$$ [X,rY]=r[X,Y]+\rho(X)(r)Y $$
for any $X,Y\in L$ and $r\in R$. We shall write $\rho(X)(r)=X(r)$.

Note that in \cite{Hue} and \cite{Rin},
a Lie-Rinehart algebra over a fixed $R$ is
referred to as $(k,R)$-Lie algebra.

\begin{ex} \label{ex1.1} \rm
(1)
If $R=k$, then $\Der_{k}(R)=\{0\}$, and a Lie-Rinehart
algebra over $R$ is simply a Lie algebra over $k$.

(2)
For arbitrary $R$, the Lie algebra $\Der_{k}(R)$ is
itself a Lie-Rinehart algebra over $R$ if one takes $\rho$
to be the identity map.

(3)
Let $k=\RR$. If $\pi\!:E\to M$ is a vector bundle
over a smooth manifold $M$, then the structure
of a Lie-Rinehart algebra on the $\eC(M)$-module
$\Gamma(E)$ of sections of $E$ is the same as the structure
of a Lie algebroid on $E$.

(4)
Let $L$ be a Lie-Rinehart algebra over $R$ and
$\tau\!:R\to S$ a homomorphism of unital commutative
$k$-algebras. Then we can form a Lie-Rinehart
algebra $\tau_{!}(L)$ over $S$ as the kernel
of the map $\varphi$,
$$
\xymatrix{
0 \ar[r] & \tau_{!}(L) \ar[r] &
(S\ten_{R}L)\oplus \Der_{k}(S) \ar[r]^-{\varphi} &
\Hom_{k}(R,S)\;,
}
$$
defined by
$$ \varphi(s\ten X,D)(r)=s\tau(X(r))-D(\tau(r))\;.$$
The bracket on $\tau_{!}(L)$ is given by
$$ [(s\ten X,D),(t\ten Y,E)]=
   (st\ten [X,Y] + D(t)\ten Y - E(s)\ten X , [D,E] )\;,$$
while the representation of $\tau_{!}(L)$ on $S$
is given by the projection.
In the special case where $\tau$ is a localization
$R\to R_{\fp}$ of $R$ at a prime ideal $\fp$ of $R$,
one can check that $\tau_{!}(L)$ is isomorphic
to the localization $R_{\fp}\ten_{R}L=L_{\fp}$, with the bracket
$$ [s^{-1}X,t^{-1}Y] =
   (st)^{-1}[X,Y] - s^{-1}t^{-2}X(t)Y + t^{-1}s^{-2}Y(s)X $$
and representation $\rho_{\fp}\!:L_{\fp}\to\Der_{k}(R_{\fp})$ induced
by $\rho\!:L\to\Der_{k}(R)$ and the canonical map
$R_{\fp}\ten_{R}\Der_{k}(R)\to\Der_{k}(R_{\fp})$.
This agrees with the definition given in \cite{Rin}.
\end{ex}

The Lie-Rinehart algebras over $R$ form a category
$$ \LRA_{R} \;,$$
where a morphism $\phi\!:L\to L'$ is a homomorphism
of Lie algebras over $k$ as well as a map
of $R$-modules which intertwines the representations,
$\rho'\com \phi=\rho$.
The operation $\tau_{!}$, induced by
a homomorphism $\tau\!:R\to S$
of unital commutative $k$-algebras,
is a functor $\LRA_{R}\to\LRA_{S}$.
Moreover, for another homomorphism
$\sigma\!:S\to T$ of unital commutative
$k$-algebras there is a canonical map
$\sigma_{!}(\tau_{!}(L))\to (\sigma\com\tau)_{!}(L)$.
Using these canonical maps, the
categories $\LRA_{R}$,
for varying $k$-algebras $R$, can be
assembled into one big (fibered) category
$$ \LRA\;,$$
in which a map $(R,L)\to (S,K)$
is a pair $(\tau,\phi)$, consisting of a homomorphism
of unital $k$-algebras $\tau\!:R\to S$ and a homomorphism
$\phi\!:\tau_{!}(L)\to K$ of Lie-Rinehart algebras over $S$.

\section{The universal enveloping algebra} \label{sec2}

Let $L$ be a Lie-Rinehart algebra over $R$.
The left $R$-module $R\oplus L$ has a natural
Lie algebra structure, given by
$$[(r,X), (s,Y)]=(X(s)-Y(r),[X,Y]) $$
for any $r,s\in R$ and $X,Y\in L$.
Let $U(R\oplus L)$ be its universal enveloping algebra over $k$,
obtained as the quotient of the tensor algebra
of $R\oplus L$ over $k$ with respect
to the usual ideal.
Write $i\!:R\oplus L\to U(R\oplus L)$ for the canonical inclusion
and $\bar{U}(R\oplus L)$ for
the subalgebra of $U(R\oplus L)$ generated by
$i(R\oplus L)$.
The {\em universal enveloping algebra} of the Lie-Rinehart algebra
$L$ over $R$ (see \cite{Rin})
is the quotient algebra
$$ \Univ(R,L)=\bar{U}(R\oplus L)/I $$
over $k$,
where $I$ is the two-sided
ideal in $\bar{U}(R\oplus L)$
generated by the elements
$i(s,0)\cdot i(r,X)-i(sr,sX)$,
for all $r,s\in R$ and $X\in L$. The natural map
$\iota_{R}\!:R\to\Univ(R,L)$, $r\mapsto i(r,0)+I$,
is a homomorphism of unital $k$-algebras, while
$\iota_{L}\!:L\to\Univ(R,L)$,
$X\mapsto i(0,X)+I$, is a homomorphism of Lie algebras.
Furthermore, we have $\iota_{R}(r)\iota_{L}(X)=\iota_{L}(rX)$
and $[\iota_{L}(X),\iota_{R}(r)]=\iota_{R}(X(r))$.

The universal enveloping algebra $\Univ(R,L)$ is characterized by the
following universal property:
if $A$ is any unital $k$-algebra,
$\kappa_{R}\!:R\to A$ a homomorphism of unital $k$-algebras
and $\kappa_{L}\!:L\to A$ a homomorphism of Lie algebras
such that $\kappa_{R}(r)\kappa_{L}(X)=\kappa_{L}(rX)$
and $[\kappa_{L}(X),\kappa_{R}(r)]=\kappa_{R}(X(r))$
for any $r\in R$ and $X\in L$, then
there exists a unique homomorphism of
unital algebras $f\!:\Univ(R,L)\to A$ such that
$f\com \iota_{R}=\kappa_{R}$ and
$f\com \iota_{L}=\kappa_{L}$.

In particular, the universal property of $\Univ(R,L)$
implies that there exists a unique
representation
$$ \varrho\!:\Univ(R,L)\to \End_{k}(R) $$
such that $\varrho\com\iota_{L}=\rho$ and
$\varrho\com\iota_{R}$ is the canonical
representation given by the
multiplication in $R$. Since
the canonical
representation of $R$ is faithful, we see that
the map $\iota_{R}$ is injective.
We shall therefore identify $\iota_{R}(R)$ with $R$,
$\iota_{R}(r)=r$.
We shall often denote
$\iota_{L}(X)$ by $\iota(X)$ or simply by $X$.
In this notation, the algebra $\Univ(R,L)$
is generated by elements $X\in L$ and $r\in R$,
while $r\cdot X=rX$ and
$[X,r]=X\cdot r - r\cdot X=X(r)$ in $\Univ(R,L)$.
As a $k$-linear space, $\Univ(R,L)$ is generated by $R$ and
the powers $\iota(L)^{n}$, $n=1,2,\ldots\,$.
The algebra $\Univ(R,L)$ also has a natural filtration
$$ R=\Univ_{(0)}(R,L)\subset\Univ_{(1)}(R,L)\subset\Univ_{(2)}(R,L)\subset\cdots\;,$$
where $\Univ_{(n)}(R,L)$ is spanned by $R$ and the
powers $\iota(L)^{m}$, for $m=1,2,\ldots,n$.
We define the associated graded algebra as
$$ \gr(\Univ(R,L))=\bigoplus_{n=0}^{\infty} \Univ_{(n)}(R,L)/\Univ_{(n-1)}(R,L)\;,$$
where we take $\Univ_{(-1)}(R,L)=\{0\}$.
It is a commutative unital algebra over $R$.

\begin{ex} \rm \label{ex2.1}
(1) Let $V$ be a left $R$-module. With zero bracket and representation,
$V$ is a Lie-Rinehart algebra. The corresponding universal enveloping algebra
$\Univ(R,V)$ is in this case the symmetric algebra $S_{R}(V)$
(see the appendix below).

(2) For any Lie algebra over $k$, the universal enveloping algebra
$\Univ(k,L)$ is the classical
universal enveloping algebra $U(L)$ of $L$.

(3) If $G$ is a Lie groupoid over a smooth manifold $M$
and $\fg$ its Lie algebroid, then the algebra of right
invariant tangential differential operators on $G$ is
a concrete model for the universal enveloping algebra $\Univ(\eC(M),\Gamma(\fg))$
(see \cite{NWX}).
\end{ex}

As for the classical universal enveloping algebra of a Lie algebra,
there is a Poincar\'{e}-Birkhoff-Witt theorem for the
universal enveloping algebra of a Lie-Rine\-hart algebra \cite{Rin}:
if the Lie-Rinehart algebra $L$ is projective as a left $R$-module,
then the
natural map $\theta\!:S_{R}(L)\to\gr(\Univ(R,L))$ is an isomorphism of
algebras. In particular, this implies that
$\iota_{L}\!:L\to\Univ(R,L)$ is in this case injective.

\section{Rinehart bialgebras} \label{sec3}

The universal enveloping algebra of a Lie algebra
is a Hopf algebra, as is the group ring of a discrete group.
In this section we will identify the algebraic structure
common to the universal enveloping algebra of a Lie-Rinehart algebra
and the convolution algebra of an \'{e}tale groupoid. This
structure has occurred in the literature under various names,
see \cite{Kap,Lu,Mal,Mrcun,Mrcun2,Tak,Xu}. We suggest the name
{\em Rinehart bialgebra}.

Let $R$ be a unital commutative algebra over $k$ as before.
All modules considered will be unital left $R$-modules, and
$\ten_{R}$ will always denote the tensor product of left
$R$-modules.
Suppose that $A$ is a unital $k$-algebra
which {\em extends} $R$, i.e.\ such that
$R$ is a unital subalgebra of $A$.
In particular, $A$ is an $R$-$R$-bimodule.
The left $R$-module $A\ten_{R}A$
(tensor product of $A$, viewed as a left $R$-module, with itself)
is also a right $R$-module in two ways.
Observe, that $A\ten_{R}A$ is not necessarily an algebra in a
natural way unless $R$ lies in the centre of $A$.
Following \cite{Kap},
we denote by $A\bar{\ten}_{R}A$ the submodule of
$A\ten_{R}A$ given by the kernel of the map $\vartheta$,
$$
\xymatrix{
0 \ar[r] & A\bar{\ten}_{R}A \ar[r] &
A\ten_{R}A \ar[r]^-{\vartheta} & \Hom_{k}(R,A\ten_{R}A)\;,
}
$$
defined by $\vartheta(a\ten b)(r)=ar\ten b - a\ten br$.
The $R$-module $A\bar{\ten}_{R}A$ has a natural structure of
a $k$-algebra.
If $R$ is in the centre of $A$, then
$A\bar{\ten}_{R}A=A\ten_{R}A$.

\begin{dfn} \rm \label{dfn3.1}
A {\em Rinehart bialgebra over} $R$
is a unital $k$-algebra $A$ which extends $R$,
with a {\em compatible} structure of a cocommutative
coalgebra in the category of left $R$-modules.

If we denote the comultiplication
by $\cm\!:A\to A\ten_{R} A$ and the counit
by $\cu\!:A\to R$, the compatibility conditions
are
\begin{enumerate}
\item [(i)] $\cm(A)\subset A\bar{\ten}_{R}A$,
\item [(ii)] $\cu(1)=1$,
\item [(iii)] $\cm(1)=1\ten 1$,
\item [(iv)] $\cu(ab)=\cu(a\cu(b))$, and
\item [(v)] $\cm(ab)=\cm(a)\cm(b)$
\end{enumerate}
for any $a,b\in A$.
\end{dfn}

Observe that the condition (v) makes sense because
of (i) and the fact that $A\bar{\ten}_{R}A$ is a $k$-algebra.
Note, however, that (iv) does not express that $\cu$ is
an algebra map.

The Rinehart bialgebras over $R$ form a category
$$ \RBA_{R}\;, $$
in which a morphism $f\!:A\to B$ is a map
which is at the same time a homomorphism of $k$-algebras
with unit and a homomorphism of coalgebras in the category
of left $R$-modules.

\begin{ex} \label{ex3.2} \rm
(1)
If a unital $k$-algebra $A$ extends $R$ such that $R$ lies in the centre
of $A$, a Rinehart bialgebra structure on $A$
is the same as an ordinary $R$-bialgebra structure on $A$.

(2)
Let $\Univ(R,L)$ be the universal enveloping algebra of
a Lie-Rinehart algebra $L$ over $R$ (Section \ref{sec2}).
The universal property implies that there exists a unique
homomorphism of algebras
$$ \cm\!:\Univ(R,L)\to \Univ(R,L)\bar{\ten}_{R}\Univ(R,L)
   \subset \Univ(R,L)\ten_{R}\Univ(R,L) $$
such that $\cm(r)=1\ten r=r\ten 1$ and
$\cm(X)=1\ten X + X\ten 1$ for any $r\in R$ and $X\in L$.
With the counit $\cu\!:\Univ(R,L)\to R$ given by
$$ \cu(u)=\varrho(u)(1) \;,$$
one can check that $\Univ(R,L)$ is a Rinehart bialgebra over $R$.
Furthermore, a morphism of Lie-Rinehart algebras
$\phi\!:L\to L'$ induces, by the universal property,
a morphism of Rinehart bialgebras
$\Univ(R,\phi)\!:\Univ(R,L)\to\Univ(R,L')$, and this gives a functor
$$ \Univ\!:\LRA_{R}\to\RBA_{R} \;.$$

(3)
The convolution algebra of
smooth functions with compact support on
an \'{e}tale Lie groupoid $G$
over a compact manifold $M$
is a Rinehart bialgebra over $\eC(M)$
(see \cite{Mrcun,Mrcun2}).
\end{ex}

Let $A$ be a Rinehart bialgebra over $R$. Then $A$ splits
as an $R$-module as
$$ A=R\oplus\bar{A}\;,$$
where $\bar{A}=\ker\cu$. The $R$-submodule $\bar{A}$
is a subalgebra of $A$ and carries a
cocommutative coassociative coproduct
$\bar{\cm}\!:\bar{A}\to \bar{A}\ten_{R}\bar{A}$, defined by
$$ \bar{\cm}(a)=\cm(a)- a\ten 1 - 1\ten a\;.$$
The bialgebra $A$ can be reconstructed from $\bar{A}$, $\bar{\cm}$
and the multiplication on $\bar{A}$.

The $R$-module $\bar{A}$ has a filtration
$$ \{0\}=\bar{A}_{0} \subset \bar{A}_{1} \subset \bar{A}_{2} \subset \cdots $$
with $\bar{A}_{n}=\ker\bar{\cm}^{(n)}$, where
$\bar{\cm}^{(n)}$ denotes the iterated coproduct
$\bar{A}\to \bar{A}\ten_{R}\cdots\ten_{R}\bar{A}$ ($n+1$ copies).
We refer to this filtration as the {\em primitive filtration}
of $\bar{A}$, and also write $A_{n}=R\oplus\bar{A}_{n}$ ($n\geq 0$).
The
submodule $\bar{A}_{1}$ is called the submodule of {\em primitive}
elements, and also denoted by $\Prim(A)$.
We observe that $\Prim(A)$ is a Lie-Rinehart algebra over $R$. Its
Lie bracket is given by the commutator in $A$, and its representation
$\rho\!:\Prim(A)\to \Der_{k}(R)$ is given by
$\rho(a)(r)=\cu(a r)$.
Indeed, note that
\begin{equation*}
\begin{split}
a r      & = (\cu\ten 1)(\cm(a  r)) \\
         & = (\cu\ten 1)(\cm(a)  \cm(r)) \\
         & = (\cu\ten 1)((a\ten 1 + 1\ten a) (r\ten 1)) \\
         & = (\cu\ten 1)(a r\ten 1 + r\ten a) \\
         & = \cu(a  r)+r a\;,
\end{split}
\end{equation*}
and from this it follows easily that $\rho$
is an $R$-linear homomorphism of Lie algebras.

We call $A$ (or $\bar{A}$) {\em cocomplete} if
$A=\bigcup_{n=0}^{\infty}A_{n}$
(or $\bar{A}=\bigcup_{n=0}^{\infty}\bar{A}_{n}$).
The Rinehart algebra $A$ is {\em graded projective}
if each of the subquotients
$A_{n+1}/A_{n}=\bar{A}_{n+1}/\bar{A}_{n}$
is a projective $R$-module.
We will use these notions in
several theorems stated below.

\begin{ex} \label{ex3.3} \rm
It follows from the Poincar\'{e}-Birkhoff-Witt theorem
that the primitive filtration of the
universal enveloping algebra $\Univ(R,L)$, associated to
a Lie-Rinehart algebra $L$ over $R$,
coincide with its natural filtration
if $L$ is projective as a left $R$-module.
Furthermore,
the universal enveloping algebra $\Univ(R,L)$ is in this case
cocomplete and graded projective.
\end{ex}

\begin{rem} \rm \label{rem3.4}
Let $A$ be a cocomplete graded projective
Rinehart bialgebra over $R$. In particular, this implies that
$A$ and all $A_{n}$ are projective $R$-modules.
We write $\gr(\bar{A})=\bigoplus_{n=1}^{\infty}\gr_{n}(\bar{A})$,
where $\gr_{n}(\bar{A})=\bar{A}_{n}/\bar{A}_{n-1}$.
There is a cocommutative coassociative comultiplication
$\bar{\cm}^{\mathrm{gr}}$ on $\gr(\bar{A})$,
induced by $\bar{\cm}$, such that
$$ \bar{\cm}^{\mathrm{gr}}(\gr_{n}(\bar{A}))\subset
   \bigoplus_{p+q=n}\gr_{p}(\bar{A})\ten \gr_{q}(\bar{A}) $$
and $\ker(\bar{\cm}^{\mathrm{gr}})=\gr_{1}(\bar{A})$.
Furthermore, the non-counital coalgebra
$\gr(\bar{A})$ is cocomplete, i.e.\
$\bigcup_{n=1}^{\infty}\ker((\bar{\cm}^{\mathrm{gr}})^{(n)})=\gr(\bar{A})$
(see the appendix).
Note that
any morphism $f\!:A\to B$ of cocomplete graded projective
Rinehart bialgebras over $R$ induces a
morphism of non-counital coalgebras
$\gr(\bar{f})\!:\gr(\bar{A})\to\gr(\bar{B})$.
\end{rem}

\section{A Cartier-Milnor-Moore theorem} \label{sec4}

We have already seen that the universal enveloping algebra construction
defines a functor
$$ \Univ\!:\LRA_{R}\to\RBA_{R}\;.$$
In the other direction, there is a functor
$$ \Prim\!:\RBA_{R}\to\LRA_{R}\;,$$
which assigns to a Rinehart bialgebra $A$
its Lie-Rinehart algebra $\Prim(A)$ of
primitive elements.

\begin{theo} \label{theo4.1}
The functor $\Univ$ is left adjoint to $\Prim$.
Furthermore, the functors $\Univ$ and $\Prim$ restrict
to an equivalence between the full subcategory
of Lie-Rinehart algebras over $R$
which are projective as left $R$-modules
and that
of cocomplete graded projective Rinehart bialgebras
over $R$.
\end{theo}

The second part of the theorem in particular implies
the following property of the counit of the adjunction,
which is an analogue of the Cartier-Milnor-Moore
theorem for Hopf algebras:

\begin{cor} \label{cor4.2}
Let $A$ be a Rinehart bialgebra over $R$.
If $A$ is cocomplete and graded projective, then
there is a canonical isomorphism of
Rinehart bialgebras
$\Univ(R,\Prim(A))\to A$.
\end{cor}

\begin{proof}[Proof of Theorem \ref{theo4.1}]
For a Lie-Rinehart algebra $L$ over $R$, the
canonical map $L\to\Univ(R,L)$ clearly lands
in the submodule of primitive elements, and this
defines the unit of the adjunction,
$$ \alpha_{L}\!:L\to\Prim(\Univ(R,L))\;.$$
For a Rinehart bialgebra $A$ over $R$, the inclusion
$\Prim(A)\to A$ induces, by the universal property
of the universal enveloping algebra, a canonical algebra map
$$ \beta_{A}\!:\Univ(R,\Prim(A))\to A\;,$$
which in fact is clearly a map of Rinehart
bialgebras. This defines the counit of the adjunction.

The first part of the theorem states that these two maps
satisfy the triangular identities
$$ \Prim(\beta_{A})\com\alpha_{\Prim(A)}=\id_{\Prim(A)} $$
and
$$ \beta_{\Univ(R,L)}\com\Univ(R,\alpha_{L})=\id_{\Univ(R,L)}\;. $$
These both hold by the (uniqueness part of the) universal
property of the universal enveloping algebra.

If $L$ is projective as a left $R$-module,
then the Poincar\'{e}-Birkhoff-Witt theorem
implies that $\Univ(R,L)$ is graded projective and cocomplete
(because the natural and primitive filtrations coincide),
and that $\alpha_{L}$ is an isomorphism.
For a graded projective cocomplete Rinehart bialgebra
$A$ over $R$, the $R$-module $\Prim(A)$ is obviously
projective, and it remains to show that in this case the map
$\beta=\beta_{A}$ is an isomorphism.

It suffices to prove that the map of reduced
non-counital coalgebras
$$ \bar{\beta}\!:\bar{\Univ}(R,\Prim(A))\to\bar{A} $$
is an isomorphism.
For this, in turn, it is enough to show that the induced map
of non-counital coalgebras
$$ \gr(\bar{\beta})\!:\gr(\bar{\Univ}(R,\Prim(A)))\to\gr(\bar{A}) $$
is an isomorphism.
The projection
$\gr(\bar{A})\to \gr_{1}(\bar{A})=\Prim(A)$ induces a map
of non-counital coalgebras $\gamma\!:\gr(\bar{A})\to \bar{S}_{R}(\Prim(A))$, by
the universal property of $\bar{S}_{R}(\Prim(A))$ (Proposition \ref{propA.2}).
Now consider the diagram
$$
\xymatrix{
\gr(\bar{\Univ}(R,\Prim(A))) \ar[r]^-{\gr(\bar{\beta})} \ar[dr]_-{\bar{\theta}^{-1}} &
\gr(\bar{A}) \ar[d]^{\gamma} \\
& \bar{S}_{R}(\Prim(A))
}
$$
where $\bar{\theta}^{-1}$ is the inverse of the Poincar\'{e}-Birkhoff-Witt
isomorphism $\theta$ (Section \ref{sec2}) restricted to
$\bar{S}_{R}(\Prim(A))$.
All maps in the diagram
are maps of non-counital coalgebras.
Thus, to see that the diagram commutes, it suffices
(again by the universal property of $\bar{S}_{R}(\Prim(A))$ stated
in Proposition \ref{propA.2}) that
$$ \pr_{1}\com\gamma\com\gr(\bar{\beta})=\pr_{1}\com\bar{\theta}^{-1} $$
for the projection $\pr_{1}\!:\bar{S}_{R}(\Prim(A))\to\Prim(A)$,
which is clear from the explicit definitions.
To finish the proof,
recall from Remark \ref{rem3.4}
that $\gr(\bar{A})$ is cocomplete and $\ker(\bar{\cm}^{\mathrm{gr}})=\Prim(A)$,
so that
by Lemma \ref{lemA.1} the map
$\gamma$ is injective, while by the commutativity
of the diagram it is also surjective. Thus $\gamma$ is
an isomorphism, and hence so is $\gr(\bar{\beta})$.
\end{proof}

\appendix

\section{Cocommutative non-counital coalgebras} \label{secA}

The main goal of this appendix is to
prove some elementary properties of cocommutative non-counital
coalgebras over a ring, in particular the universal property
of the symmetric coalgebra. Although these results are
well-known, they are usually stated in the literature in the context
of coalgebras over a field (see e.g.\ \cite{Car,MilMoo,Qui,Swe}).

As before, $k$ denotes a field of characteristic $0$
and $R$ a unital commutative $k$-algebra.
A cocommutative non-counital coalgebra over $R$
is an $R$-module $C$, together with a cocommutative
coassociative comultiplication.
$\delta\!:C\to C\ten_{R}C$. (Note that
we do not assume that $C$ has a counit.)
With the obvious maps, these cocommutative non-counital
coalgebras over $R$ form a category.
For example, for a Rinehart bialgebra $A$ over $R$,
the pair $(\bar{A},\bar{\cm})$ is an object of this
category.

As discussed in this special case already, any
cocommutative non-counital coalgebra $(C,\delta)$ carries
a {\em primitive} filtration
$$ \{0\}=C_{0}\subset C_{1}\subset C_{2}\subset \cdots\;,$$
where $C_{n}$ is the kernel of the iterated coproduct
$\delta^{(n)}\!:C\to C\ten_{R}\cdots\ten_{R}C$ ($n+1$ copies).
We say that $C$ is {\em cocomplete} if $C=\bigcup_{n=0}^{\infty}C_{n}$.
If $C$ is cocomplete and $C_{n}$ are all flat $R$-modules,
then $C$ is a flat $R$-module as well.
Notice that in this case
$\delta\!:C\to C\ten_{R}C$ maps each
$C_{n+1}$ into $C_{n}\ten_{R}C_{n}$.

\begin{lem} \label{lemA.1}
Let $f\!:C\to D$ be a morphism of
cocomplete cocommutative non-counital
coalgebras over $R$, and assume that all
the submodules $C_{n}$, $D_{n}$ are flat.
If $f_{1}=f|_{C_{1}}\!:C_{1}\to D_{1}$
is injective, then so is $f\!:C\to D$.
\end{lem}

\begin{proof}
We prove by induction that $f_{n}=f|_{C_{n}}\!:C_{n}\to D_{n}$
is injective. Assuming this has been proved for $f_{n}$,
injectivity of $f_{n+1}$ follows from the diagram
$$
\xymatrix{
C_{n+1} \ar[d]_-{\delta_{C}} \ar[r]^{f_{n+1}} &
D_{n+1} \ar[d]^-{\delta_{D}} \\
C_{n}\ten_{R}C_{n} \ar[r]^-{f_{n}\ten f_{n}} &
D_{n}\ten_{R}D_{n}
}
$$
since the map at the bottom is injective by
the flatness assumption.
Indeed, we have $\ker(f_{n+1})\subset\ker(\delta_{C})=C_{1}$
and $C_{1}\cap\ker(f)=\{0\}$ by assumption.
\end{proof}

Let $V$ be an $R$-module. We denote by
$S_{R}(V)$ the symmetric algebra on $V$. It is a graded algebra,
where $S^{n}_{R}(V)$ is the space
$(V\ten_{R}\cdots\ten_{R}V)_{\Sigma_{n}}$
of coinvariants of the $n$-fold
tensor product. The algebra $S_{R}(V)$ is the free commutative
unital $R$-algebra on $V$. Therefore, the map
$\cm\!:V\to S_{R}(V)\ten_{R}S_{R}(V)$, given by
$\cm(v)=1\ten v + v\ten 1$, extends uniquely to a unital algebra
homomorphism
$\cm\!:S_{R}(V)\to S_{R}(V)\ten_{R}S_{R}(V)$, giving
$S_{R}(V)$ the structure of a commutative and cocommutative bialgebra
over $R$.
(In fact, this is the Rinehart bialgebra $\Univ(R,L)$, where
$L$ is $V$ viewed as a Lie-Rinehart algebra with
zero bracket and representation.)

The bialgebra $S_{R}(V)$ also has a universal property as a coalgebra,
most easily stated in terms of the $R$-module
$\bar{S}_{R}(V)=\bigoplus_{n=1}^{\infty}S^{n}_{R}(V)$, which
is a cocomplete cocommutative non-counital
$R$-algebra with coproduct $\bar{\cm}(c)=\cm(c)-1\ten c- c\ten 1$.

\begin{prop} \label{propA.2}
Let $(C,\delta)$ be a cocomplete cocommutative non-counital
coalgebra over $R$ and $V$ an $R$-module. Then any
morphism of $R$-modules $C\to V$ is the first component of
a unique homomorphism $C\to \bar{S}_{R}(V)=\bigoplus_{n=1}^{\infty}S^{n}_{R}(V)$
of non-counital coalgebras over $R$.
\end{prop}

\begin{proof}
Let $g_{1}\!:C\to V=S^{1}_{R}(V)$ be an $R$-linear map.
We will define a map $g_{n}\!:C\to S^{n}_{R}(V)$, for each
$n\geq 2$, in such a way that these together form a coalgebra
map $g=(g_{n})\!:C\to \bar{S}_{R}(V)=\bigoplus_{n=1}^{\infty}S^{n}_{R}(V)$.
Note that for $g$ to be a coalgebra map, there can be at most one
such $g_{n}$, since it will have to make the diagram
$$
\xymatrix{
C \ar[d]_-{\delta^{(n+1)}} \ar[r]^-{g_{n}} &
S^{n}_{R}(V) \ar[d]^-{\bar{\cm}^{(n+1)}} \\
C^{\ten_{R}n} \ar[r]^-{g_{1}^{\ten n}} &
{S^{1}_{R}(V)}^{\ten_{R}n}
}
$$
commute. In fact, we can use this diagram to define $g_{n}$,
because $C$ is cocommutative and $\bar{\cm}^{(n+1)}$ gives
an isomorphism between $S^{n}_{R}(V)$ and the subspace
of $\Sigma_{n}$-invariants
of ${S^{1}_{R}(V)}^{\ten_{R}n}$.
(Here we used the assumption that the characteristic of $k$ is $0$.)

The cocompleteness of $C$ implies that $g$ is well defined.
To see that $g$ is indeed a map of coalgebras,
we show by induction on $n$ that the diagram
\begin{equation} \label{diagA.3}
\begin{split}
\xymatrix{
C \ar[d]_-{\delta} \ar[rr]^-{g_{n}} & &
S^{n}_{R}(V) \ar[d]^-{\bar{\cm}} \\
C\ten_{R}C \ar[rr]^-{(g_{p}\ten g_{n-p})} & &
\bigoplus_{0<p<n} (S^{p}_{R}(V)\ten S^{n-p}_{R}(V))
}
\end{split}
\end{equation}
commutes.
For $n=1$ there is nothing to prove, and for $n=2$
this holds by definition of $g_{2}$. Take $n>2$, and suppose that
the commutativity of (\ref{diagA.3}) has been proved for all
$m<n$. Now consider the following diagram
\begin{equation} \label{diagA.4}
\begin{split}
\xymatrix{
C \ar[rrr]^-{g_{n}} \ar[d]_-{\delta} &&&
S^{n}_{R}(V) \ar[d]^-{\bar{\cm}} \\
C^{\ten_{R} 2} \ar[rrr]^-{(g_{p_{1}}\ten g_{p_{2}})}
\ar@<-2pt>[d]_-{\delta\ten 1} \ar@<2pt>[d]^-{1\ten\delta} &&&
\bigoplus_{p_{1}+p_{2}=n} (S^{p_{1}}\ten_{R} S^{p_{2}})
\ar@<-2pt>[d]_-{\bar{\cm}\ten 1} \ar@<2pt>[d]^-{1\ten\bar{\cm}} \\
C^{\ten_{R} 3} \ar[rrr]^-{(g_{p_{1}}\ten g_{p_{2}}\ten g_{p_{3}})}
\ar@<-4pt>[d]_-{\delta\ten 1\ten 1}
\ar[d]
\ar@<4pt>[d]^-{1\ten 1\ten\delta}
&&&
\bigoplus_{p_{1}+p_{2}+p_{3}=n} (S^{p_{1}}\ten_{R} S^{p_{2}}\ten_{R} S^{p_{3}})
\ar@<-4pt>[d]_-{\bar{\cm}\ten 1\ten 1}
\ar[d]
\ar@<4pt>[d]^-{1\ten 1\ten\bar{\cm}} \\
\vdots \ar@{}[rrr]_-{.\;\;.\;\;.}
\ar@<-8pt>[d]
\ar@<-4pt>[d]^-{\cdots}
\ar@<8pt>[d]
&&& \vdots
\ar@<-8pt>[d]
\ar@<-4pt>[d]^-{\cdots}
\ar@<8pt>[d] \\
C^{\ten_{R}n} \ar[rrr]^-{g_{1}^{\ten n}} &&&
(S^{1}_{R}(V))^{\ten_{R}n}
}
\end{split}
\end{equation}
where all $p_{i}\geq 1$ and $S^{p_{i}}=S^{p_{i}}_{R}(V)$.
Our aim is to prove that the top square in
(\ref{diagA.4}) commutes. By the induction hypothesis,
the two parallel squares at the second level commute, as
do the parallel squares at each lower level.
Furthermore, the outer square commutes by the definition
of $g_{n}$, as the vertical compositions $C\to C^{\ten_{R}n}$
all agree and define $\delta^{(n+1)}$, and
similarly for all the vertical compositions
$S^{n}_{R}(V)\to (S^{1}_{R}(V))^{\ten_{R}n}$.
Thus, to show that the top square commutes, it suffices
to prove that the joint kernel of all possible
vertical compositions
$$ \bigoplus_{p_{1}+p_{2}=n} (S^{p_{1}}\ten_{R} S^{p_{2}})
   \to (S^{1}_{R}(V))^{\ten_{R}n} $$
is zero, which can easily be checked
(note that $\bar{\cm}\!:S^{n}_{R}(V)\to
\bigoplus_{p_{1}+p_{2}=n} (S^{p_{1}}\ten_{R} S^{p_{2}})$
is split injective map of left $R$-modules if $n\geq 2$).
\end{proof}

\end{document}